\documentclass[12pt]{article}

\usepackage{graphicx,tikz,color,url}
\usepackage{amsmath,amssymb,latexsym}

\newtheorem{theorem}{Theorem}[section]

\newtheorem{lemma}[theorem]{Lemma}
\newtheorem{proposition}[theorem]{Proposition}

\newtheorem{claim}{Claim}

\newtheorem{example}{Example}

\newcommand{\proof}{\noindent{\bf Proof.\ }}
\newcommand{\qed}{\hfill $\square$ \bigskip}
\newcommand{\smallqed}{\hfill {\tiny $\square$}\bigskip}
\newcommand{\pch}{\chi_{\rho}}
\newcommand{\Spch}{\chi_{S}}
\newcommand{\diam}{{\rm diam}}

\begin{document}

\title{$S$-packing chromatic vertex-critical graphs}

\author{P\v{r}emysl Holub$^{a}$\thanks{Email: \texttt{holubpre@kma.zcu.cz}} 
\and Marko Jakovac$^{b,c}$\thanks{Email: \texttt{marko.jakovac@um.si}}
\and Sandi Klav\v zar $^{d,b,c}$\thanks{Email: \texttt{sandi.klavzar@fmf.uni-lj.si}}
}
\maketitle

\begin{center}
$^a$ Faculty of Applied Sciences, University of West Bohemia, Pilsen, Czech Republic\\
\medskip

$^b$ Faculty of Natural Sciences and Mathematics, University of Maribor, Slovenia\\
\medskip

$^c$ Institute of Mathematics, Physics and Mechanics, Ljubljana, Slovenia\\
\medskip

$^d$ Faculty of Mathematics and Physics, University of Ljubljana, Slovenia\\
\end{center}

\begin{abstract}
For a non-decreasing sequence of positive integers  $S = (s_1,s_2,\ldots)$, the {\em $S$-packing chromatic number} $\chi_S(G)$ of $G$ is the smallest integer $k$ such that the vertex set of $G$ can be partitioned into sets $X_i$, $i \in [k]$, where vertices in $X_i$ are pairwise at distance greater than $s_i$. In this paper we introduce $S$-packing chromatic vertex-critical graphs, $\chi_{S}$-critical for short, as the graphs in which $\chi_{S}(G-u)<\chi_{S}(G)$ for every $u\in V(G)$. This extends the earlier concept of the packing chromatic vertex-critical graphs. We show that if $G$ is $\chi_{S}$-critical, then the set $\{ \chi_{S}(G)-\chi_{S}(G-u); \, u\in V(G) \}$ can be almost arbitrary. If $G$ is $\chi_{S}$-critical and $\chi_{S}(G)=k$ ($k\in \mathbb{N}$), then $G$ is called $k$-$\chi_{S}$-critical. We characterize $3$-$\chi_{S}$-critical graphs and partially characterize $4$-$\chi_{S}$-critical graphs when $s_1>1$. We also deal with $k$-$\chi_{S}$-criticality of trees and caterpillars.
\end{abstract}


\medskip

\noindent
{\bf AMS Subj.\ Class.\ (2010)}: 05C57, 05C69

%%%%%%%%%%%%%%%%%%%%%%%%%%%%%%%%%%%%%%%%%%%%%%%%%%%%%
%%%%%%%%%%%%%%%%%%%%%%%%%%%%%%%%%%%%%%%%%%%%%%%%%%%%%
\section{Introduction}
\label{sec:intro}
%%%%%%%%%%%%%%%%%%%%%%%%%%%%%%%%%%%%%%%%%%%%%%%%%%%%%
%%%%%%%%%%%%%%%%%%%%%%%%%%%%%%%%%%%%%%%%%%%%%%%%%%%%%

The {\em packing chromatic number} $\pch(G)$ of a graph $G$ is the smaller integer $k$ for which there exists a mapping $c:V(G)\rightarrow [k] = \{1,\ldots, k\}$, such that if $c(u) = c(v) = \ell$ for $u\ne v$, then $d_G(u,v) > \ell$. (Here and later $d_G(u,v)$ denotes the shortest-path distance between $u$ and $v$ in $G$.) Such a map $c$ is called a {\em packing $k$-coloring}. This concept was introduced in~\cite{goddard-2008}, named with the present names in~\cite{bresar-2007}, and extensively studied afterwards. See~\cite{balogh-2018, balogh-2019, barnaby-2017, bresar-2017, ght-2019, korze-2019, shao-2015}, references therein, as well as~\cite{czap-2019} for a variant of a facial packing coloring. 

A far reaching generalization of the packing chromatic number, formally introduced by Goddard and Xu in~\cite{goddard-2012}, but being implicitly present already in~\cite{goddard-2008}, is the following. Let $S = (s_1,s_2,\ldots)$ be a non-decreasing sequence of positive integers. An {\em $S$-packing $k$-coloring} of a graph $G$ is a mapping $c:V(G)\rightarrow [k]$, such that if $c(u) = c(v) = \ell$ for $u\ne v$, then $d_G(u,v) > s_\ell$. The {\em $S$-packing chromatic number} $\chi_S(G)$ of $G$ is the smallest integer $k$ such that $G$ admits an  $S$-packing $k$-coloring. Note that if $S = (1,1,1,\ldots)$, then we are talking about the standard proper vertex coloring and if $S=(1,2,3,\ldots)$, then we deal with the packing coloring. For investigations of $S$-packing colorings see~\cite{gastineau-2015a, gastineau-2015b, gastineau-2016, gastineau-2019, goddard-2014, maarouf-2017}. 

Now, in~\cite{klavzar-2019} the {\em packing chromatic vertex-critical graphs} were introduced as the graphs $G$ for which $\pch(G-u) < \pch(G)$ holds for every $u\in V(G)$. In this paper we are interested if (and how) the results from~\cite{klavzar-2019} extend from packing colorings to $S$-packing colorings. For this sake we say that $G$ is {\em $S$-packing chromatic vertex-critical} if $\Spch(G-u) < \Spch(G)$ holds for every $u\in V(G)$. 

We proceed as follows. In Section \ref{sec:prelim}, we list some known results and prove two statements which will be used in the rest of the paper. In the subsequent section we consider the effect of vertex removal on the $S$-packing chromatic number and prove two realization theorems. Setting $\Delta_{\Spch}(G) = \{\Spch(G) - \Spch(G-u):\ u\in V(G)\}$, the first of these results asserts that if  $S = (1^\ell,2^{\infty})$, $\ell \ge 1$, and $A =\{1,a_1,\ldots,a_k\}$, $k \ge 1$, is a set of positive integers, then there exists a $\Spch$-critical graph $G$ such that $\Delta_{\Spch}(G) = A$. In Section \ref{sec:3-critical}, we give a complete characterization of $3$-$\Spch$-critical graphs while in Section \ref{sec:4-critical} we partially characterize $4$-$\Spch$-critical graphs for packing sequences with $s_1>1$. Finally, in Section \ref{sec:tress} we show that a $k$-$\Spch$-critical tree exists for any positive integer $k$, investigate $k$-$\Spch$-criticality of caterpillars for sequences $S=(1,s_2^{k-1})$, and give some examples of such critical cattepillars.

%%%%%%%%%%%%%%%%%%%%%%%%%%%%%%%%%%%%%%%%%%%%%%%%%%%%%
%%%%%%%%%%%%%%%%%%%%%%%%%%%%%%%%%%%%%%%%%%%%%%%%%%%%%
\section{Preliminaries}
\label{sec:prelim}
%%%%%%%%%%%%%%%%%%%%%%%%%%%%%%%%%%%%%%%%%%%%%%%%%%%%%
%%%%%%%%%%%%%%%%%%%%%%%%%%%%%%%%%%%%%%%%%%%%%%%%%%%%%

The order of a graph $G$ will be denoted with $n(G)$. A graph consisting of a triangle and an edge with one end vertex on the triangle will be denoted by $Z_1$ (see also the top-right graph in Fig.~\ref{fig:4-crit}).

We will be interested in non-decreasing finite or infinite sequences $S = (s_1,s_2,\ldots)$ of positive integers, but exclude the constant sequence $(1,1,\ldots)$ because it leads to the chromatic number for which critical graphs are already well-studied, cf.~\cite{jensen-1995}. For any other sequence $S$ we will say that $S$ is a {\em packing sequences}. If in a packing sequence a term $i$ is repeated $\ell$ times, we will abbreviate the corresponding subsequence by $i^\ell$. For instance, if $S = (1,\ldots, 1, s_{\ell + 1}, \ldots)$, that is, if $S$ starts with $\ell$ terms equal to $1$, then we will write $S = (1^\ell, s_{\ell + 1}, \ldots)$.  Moreover, we will use the same convention for infinite constant subsequences. For example, $(1^\ell, 2, 2, \ldots)$ will be abbreviated $(1^\ell, 2^{\infty})$.

Let $\alpha_k(G)$, $k\ge 1$, denote the maximum number of vertices of a graph $G$ that can be properly colored using $k$ colors, that is, the cardinality of a largest $k$-independent set of $G$. We now recall a series of results from~\cite{goddard-2012} that will be used later.

\begin{lemma} {\rm (\cite[Observation 2]{goddard-2012})}
\label{lem:subgraph}
Let $S$ be a packing sequence. If $H$ is a subgraph of $G$, then $\Spch(H) \leq \Spch(G)$.
\end{lemma}

\begin{proposition} {\rm (\cite[Proposition 4]{goddard-2012})}
\label{prp:2-col}
Let $S=(s_1,s_2,\dots)$ be a packing sequence and let $G$ be a nonempty connected graph. 
\begin{enumerate}
\item If $s_1=s_2=1$, then $\Spch(G)=2$ if and only if $G$ is bipartite.
\item If $s_1=1$ and $s_2>1$, then $\Spch(G)=2$ if and only if $G$ is a star. 
\item If $s_1\geq 2$, then $\Spch(G)=2$ if and only if $G\simeq K_2$.
\end{enumerate}
\end{proposition}

\begin{proposition} {\rm (\cite[Proposition 6]{goddard-2012})}
\label{prp:Wayne-diam-2}
Let $S = (1^\ell, s_{\ell+1}, \ldots)$, where $\ell\ge 1$ and $s_{\ell+1}\ge 2$, and let $G$ be a graph with $\diam(G) = 2$. Then $\Spch(G) = n(G) - \alpha_\ell(G) + \min\{\ell, \chi(G)\}$.
\end{proposition}

\begin{proposition} {\rm (\cite[Proposition 20]{goddard-2012})}
\label{prp:3-col-222}
Let $G$ be a connected graph and $S=(2,2,2)$. Then $G$ has a $\Spch$-packing coloring if and only if $G$ is a path of any length or a cycle of lenght a multiple of $3$.
\end{proposition}

\begin{proposition} {\rm (\cite[Corollary 21]{goddard-2012})}
\label{prp:3-col-223}
Let $G$ be a graph and $S=(s_1,s_2,s_3)$, where $s_1=2$ and $s_3 \ge 3$, be a packing sequence. If $G$ has a $\Spch$-coloring, then $n(G)\leq 5$.
\end{proposition}

We conclude the preliminaries with two simple observations. 

\begin{lemma} 
\label{lem:connected}
If $S$ is a packing sequence and $G$ is a $\Spch$-critical graph, then $G$ is connected.
\end{lemma}

\proof
Since $\Spch(G)=\max_i\{\Spch(G_i)\}$, where $G_i$ are components of $G$, it follows that $G$ must have only one component provided $G$ is $\Spch$-critical.
\qed

\begin{lemma}
\label{lem:leaf}
Let $S$ be a packing sequence. If $u$ is a leaf of a graph $G$, then $\Spch(G)-1 \leq \Spch(G-u) \leq \Spch(G)$.
\end{lemma}

\proof
By Lemma~\ref{lem:subgraph}, $\Spch(G-u) \leq \Spch(G)$. Suppose that $\Spch(G-u)=k$. Then using an optimal $S$-packing coloring of $G-u$, and using color $k+1$ for the vertex $u$ in $G$, we obtain that $\Spch(G) \leq k+1$. Thus, $\Spch(G) \leq \Spch(G-u) + 1$.
\qed

%%%%%%%%%%%%%%%%%%%%%%%%%%%%%%%%%%%%%%%%%%%%%%%%%%%%%
%%%%%%%%%%%%%%%%%%%%%%%%%%%%%%%%%%%%%%%%%%%%%%%%%%%%%
\section{Vertex-deleted subgraphs of $\Spch$-critical graphs}
\label{sec:vertex-delete}
%%%%%%%%%%%%%%%%%%%%%%%%%%%%%%%%%%%%%%%%%%%%%%%%%%%%%
%%%%%%%%%%%%%%%%%%%%%%%%%%%%%%%%%%%%%%%%%%%%%%%%%%%%%

In~\cite[Theorem 3.1]{klavzar-2019} it was shown that if $G$ is a $\pch$-critical graph, then the set of differences $\Delta_{\pch}(G) = \{\pch(G) - \pch(G-u):\ u\in V(G)\}$ can be almost arbitrary. Hence the condition $\pch(G-u) < \pch(G)$ for $G$ to be $\pch$-critical cannot be replaced with  $\pch(G-u) = \pch(G) -1$. We now show with a bit more involved construction than the one from~\cite{klavzar-2019} that the same phenomenon (actually even more general) holds for all sequences of the form $(1^\ell,2^{\infty})$. More precisely, setting 
$$\Delta_{\Spch}(G) = \{\Spch(G) - \Spch(G-u):\ u\in V(G)\}$$
we have the following result. 

\begin{theorem}
\label{thm:many-values}
Let $S = (1^\ell,2^{\infty})$, $\ell \ge 1$, and let $A =\{1,a_1,\ldots,a_k\}$, $k \ge 1$, be a set of positive integers. Then there exists a $\Spch$-critical graph $G$ such that $\Delta_{\Spch}(G) = A$.
\end{theorem}

\proof
We may assume without loss of generality that $2\le a_1 < \cdots < a_k$. 

Take a cycle of length $2k+1$ on vertices $x_1, \ldots, x_{2k+1}$ (in the natural order) and additional disjoint cliques $Q_i$, $i\in [2k+1]$, where $n(Q_i) = a_{\lceil i/2\rceil} + \ell -1$ for $i\in [2k]$ and $n(Q_{2k+1}) = a_{k}  + \ell -1$. Then for every $i\in [2k+1]$ and for every vertex $x$ of $Q_i$, make $x$ adjacent to each of the vertices $x_i$, $x_{i+1}$, and $x_{i+3}, x_{i+5}, \ldots, x_{i-2}$, where indices are taken modulo $2k+1$. Denote the constructed graph with $G(\ell; a_1,\ldots, a_k)$, see Fig.~\ref{fig:G(2;2,4)} for the graph $G(2;2,4)$. 

\begin{figure}[ht!]
\centering
 \includegraphics{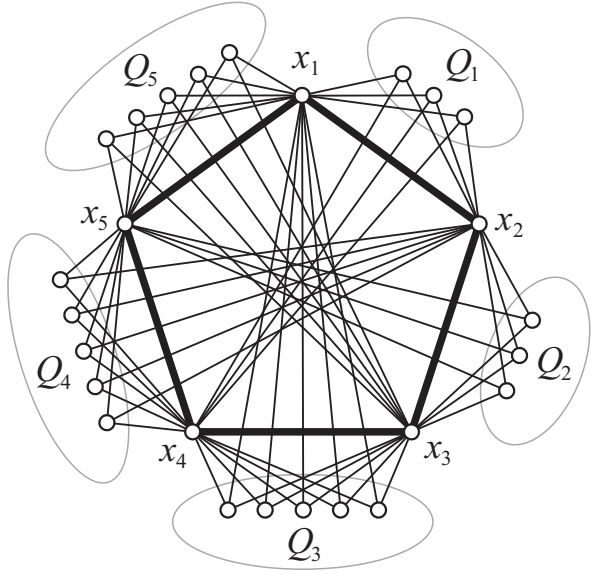}
\caption{The graph $G(2;2,4)$}
\label{fig:G(2;2,4)}
\end{figure}

Note that the vertices of the clique $Q_i$ together with the vertices $x_i$ and $x_{i+1}$ form a clique of order $a_{\lceil i/2\rceil}+ \ell + 1$. Denote this latter clique with $Q_i'$. To simplify the notation we set $G = G(\ell;a_1,\ldots, a_k)$ for the rest of the proof.

We claim first that $\diam(G) = 2$. For each $u_i\in Q_i$ and each $u_j\in Q_j$, $d_G(u_i,u_j)\leq 2$ since some of $x_j, x_{j+1}$ is adjacent to both $u_i$ and $u_j$ by the definition of $G$. Analogously, for any $x_i$ and any $u\in Q_j$, either $x_iu\in E(G)$ or $x_{i+1}u\in E(G)$ and clearly $x_ix_{i+1}\in E(G)$, implying that $d_G (x_i,u)\leq 2$. Finally, for any $x_i$ and $x_j$, either $x_ix_j\in E(G)$ or $x_i$ is adjacent to each vertex of $Q_{j-1}$ or $Q_j$, implying that $d_G(x_i,x_j)\leq 2$ since $x_ju_j\in E(G)$ for every $u_j \in Q_{j-1}\cup Q_j$. 

We have thus shown that $\diam(G) = 2$. Since clearly $\ell < \chi(G)$, Proposition~\ref{prp:Wayne-diam-2} implies that $\Spch(G) = n(G) - \alpha_\ell(G) + \ell$. Moreover, $\alpha(G) = 2k+1$ and selecting $\ell$ vertices from each of the cliques $Q_i$ we find an $\ell$-independent set of order $\ell(2k+1)$. So $\alpha_\ell(G) = \ell(2k+1)$ and consequently 
\begin{align*}
\Spch(G) & = n(G) - 2k\ell = (2k+1) + 2\sum_{i=1}^k(a_i + \ell - 1) + (a_k + \ell -1)- 2k\ell \\ 
& = 2\sum_{i=1}^ka_i + a_{k+1} + \ell\,.
\end{align*}

Let $i\in [2k+1]$ and let $x$ be an arbitrary vertex of $Q_i$. Having in mind that $n(Q_i) \ge \ell + 1$ we can repeat the above argument on the graph $G-x$ to get $\Spch(G - x) = 2\sum_{i=1}^ka_i + a_{k+1} + \ell - 1$. Hence $\Spch(G) - \Spch(G - x) = 1$.

Consider now the graph $G-x_{2i}$, where $i\in [k]$. Note that in $G$, $n(Q_{2i-1}') = n(Q_{2i}') = a_i + \ell + 1$ and that $x_{2i}$ is the unique common vertex of $Q_{2i-1}'$ and $Q_{2i}'$. If $u\in V(Q_{2i-1})$ and $v\in V(Q_{2i})$, then in $G$, the vertex $x_{2i}$ is the unique common neighbor of $u$ and $v$. It follows that $d_{G-x_{2i}}(u,v) = 3$. On the other hand, for any other pair of vertices $u'$ and $v'$ of $G-x_{2i}$ we have $d_{G-x_{2i}}(u',v') = d_{G}(u',v')$. In particular, $d_{G-x_{2i}}(x_{2i-1}, x_{2i+1}) = 2$. 

Since $\alpha_\ell(G-x_{2i}) = \ell(2k+1)$, we can select $\ell$ vertices from each of the cliques $Q_i$ to form the first $\ell$ color classes. In addition, in $V(Q_{2i-1}')-\{x_{2i}\}$ and $V(Q_{2i}')-\{x_{2i}\}$ we have $a_i-1$ pairs of vertices that are pairwise at distance $3$. We can respectively color these pairs with colors $\ell +1,\ldots, \ell + a_i - 1$. Because of the distances, every other not yet colored vertex requires its private color. Hence, with respect to the above optimal coloring of $G$, we have saved $a_i-1$ colors. Since clearly $n(G-x_{2i}) = n(G) - 1$, we thus have $\Spch(G-x_{2i}) = 2\sum_{i=1}^ka_i + a_{k+1} + \ell - (a_i-1) - 1$, which in turn implies that $\Spch(G) - \Spch(G - x_{2i}) = a_i$.

We have proved by now that $A \subseteq \Delta_{\Spch}(G)$. Finally, since $a_i < a_{i+1}$, by arguments parallel to the above arguments for $x_{2i}$ we deduce that $\Spch(G) - \Spch(G - x_{2i+1}) = a_i$. We conclude that $A = \Delta_{\Spch}(G)$.
\qed

We proceed with the case where a packing sequences contains an element which is at least $3$. To prove that the set of differences $\Delta_{\Spch}(G)$ can be almost arbitrary for this case, we can follow the same line of thought than in the proof of ~\cite[Theorem 3.1]{klavzar-2019} for packing colorings, with a few key differences in the proof. What follows is the $S$-packing coloring version of this theorem.

\begin{theorem}
\label{thm:many-values-2}
Let $S$ be a packing sequence such that there exists $\ell \ge 1$ with $s_{\ell} \ge 3$, and let $A =\{1,a_1,\ldots,a_k\}$, $k \ge 1$, be a set of positive integers. If for every $i \in [k]$ we have $\displaystyle\sum_{j=1,j \neq i}^{k}a_j \ge a_i -1$, then there exists a $\Spch$-critical graph $G$ such that $\Delta_{\Spch}(G) = A$.
\end{theorem}

\proof
Let $S = (s_1,s_2,\ldots)$ be a packing sequence and $\ell \ge 1$ the smallest positive integer such that $s_{\ell} \ge 3$.

First suppose that $k \ge 2$, and let $V(K_k)=\{x_1, \dots ,x_k\}$. We denote by $G(\ell; a_1,\ldots, a_k)$ the graph obtained from $K_k$ such that for every $i \in [k]$, a vertex of a complete graph $X_i$ of order $a_i + \ell -1$ is identified with $x_i$. Again, we simplify the notation by setting $G=G(\ell; a_1,\ldots, a_k)$ in the remainder of the proof. We first observe that 
$$n(G)=\sum_{i=1}^{k}n(X_i)=\sum_{i=1}^{k}(a_i + \ell -1)=\sum_{i=1}^{k}a_i + k(\ell -1).$$
If $c$ is a $\Spch$-coloring of $G$, then the vertices of $X_i$, $i \in [k]$, receive pairwise different colors. To be more precise, we have $|c^{-1}(j)| \leq k$ for any $j \leq \ell-1$. Moreover, since $\diam(G)=3$, $|c^{-1}(j)| \leq 1$ for any $j \ge \ell$. Since $a_i \ge 2$, and so $a_i + \ell -1 \ge \ell +1$, in each $X_i$ colors $1, \ldots, \ell -1$ can be used. Therefore, 
\begin{equation}
\label{equation1}
\Spch(G)=(\ell -1) + (n(G)- k(\ell -1))=(\ell -1) + \sum_{i=1}^{k}a_i.
\end{equation}
Since $k \ge 2$, for at least one $a_i$ we have $a_i \ge 3$. Without loss of generality we can assume that $a_1 \ge 3$. Let $u \in V(X_1)$ be an arbitrary vertex different from $x_1$. Then $G-u$ is isomorphic to $G(\ell; a_1-1, a_2, \ldots, a_k)$ (it is possible that $a_1-1=a_i$ for some $i \ge 2$). By (\ref{equation1}) we get
$$\Spch(G-u)=(\ell -1) + (a_i-1) + \sum_{i=2}^{k}a_i = \sum_{i=1}^{k}a_i + (\ell -2)=\Spch(G)-1.$$
This shows that $1 \in \Delta_{\Spch}(G)$.

Now we consider the vertex-deleted subgraph $G-x_i$, $i \in [k]$. Since $x_i$ is a cut-vertex, and using (\ref{equation1}), we have
\begin{align*}
\Spch(G-x_i)&=\max\{\Spch(K_{a_i + \ell -2}),\Spch(G(\ell; a_1, \ldots, a_{i-1}, a_{i+1}, \ldots, a_k))\}\\
									 &=\max\left\{a_i + \ell -2, (\ell -1) + \sum_{j=1,j \neq i}^{k}a_j  \right\}\\
									 &=(\ell -1) + \sum_{j=1,j \neq i}^{k}a_j,
\end{align*}
where the last inequality follows from the assumption $\displaystyle\sum_{j=1,j \neq i}^{k}a_j \ge a_i -1$. It follows that
$$\Spch(G)-\Spch(G-x_i)=\left((\ell -1) + \sum_{i=1}^{k}a_i \right) - \left((\ell -1) + \sum_{j=1,j \neq i}^{k}a_j \right)=a_i,$$
and we get $a_i \in \Delta_{\Spch}(G)$ for every $i \in [k]$.

Suppose now that $k=1$, and $A=\{1,a\}$, where $a \ge 2$. In this case, let $\ell \ge 1$ be the smallest index with $s_{\ell} \ge 2$ (not $3$ as in the previous case). Let $G$ be the graph obtained from two disjoint copies of $K_{a + \ell - 1}$ by identifying a vertex from one copy with a vertex from the other copy, and let $x$ be the identified vertex. We have $n(G)=2(a + \ell - 1)-1=2(a+\ell) - 3$ and $\Spch(G)=(\ell - 1) + (n(G) - 2(\ell - 1))=2a + \ell -2$, since $\diam(G)=2$. If $u$ is a vertex of $G$ different from $x$, then $\Spch(G-u)=(\ell - 1) + ((n(G)-1) - 2(\ell - 1))=2a + \ell -3$. Thus, $\Spch(G)-\Spch(G-u)=1$. To end the proof, we notice that $\Spch(G-x)=\Spch(K_{a+\ell-2})=a+\ell-2$, and hence $\Spch(G)-\Spch(G-x)=a$.
\qed

%%%%%%%%%%%%%%%%%%%%%%%%%%%%%%%%%%%%%%%%%%%%%%%%%%%%%
%%%%%%%%%%%%%%%%%%%%%%%%%%%%%%%%%%%%%%%%%%%%%%%%%%%%%
\section{$3$-$\Spch$-critical graphs}
\label{sec:3-critical}
%%%%%%%%%%%%%%%%%%%%%%%%%%%%%%%%%%%%%%%%%%%%%%%%%%%%%
%%%%%%%%%%%%%%%%%%%%%%%%%%%%%%%%%%%%%%%%%%%%%%%%%%%%%

If $S$ is an arbitrary packing sequence, then it is clear that $K_2$ is the unique $2$-$\Spch$-critical graph. In the following theorem we give a complete list of all $3$-$\Spch$-critical graphs with respect to a given packing sequence $S$.

\begin{theorem}
\label{thm:3-critical}
Let $S$ be a packing sequence and let $G$ be a graph.
\begin{enumerate}
\item If $S=(1,1,\ldots)$, then $G$ is $3$-$\Spch$-critical if and only if $G \in \{C_{2k+1} \, : \, k \geq 1\}$.
\item If $S=(1,s_2, \ldots)$, $s_2 \geq 2$, then $G$ is $3$-$\Spch$-critical if and only if $G\in \{C_3, C_4, P_4\}$.
\item If $S=(s_1,s_2, \ldots)$, $s_1 \geq 2$, then $G$ is $3$-$\Spch$-critical if and only if $G \in \{C_3, P_3\}$.
\end{enumerate}
\end{theorem}

\proof
Let $G$ be a $3$-$\Spch$-critical graph. Then $G$ is connected by Lemma~\ref{lem:connected}. Clearly, $n(G) \geq 3$. If $n(G) = 3$, then $G$ is either  $C_3$ or $P_3$. Clearly, $C_3$ is $3$-$\Spch$-critical for every packing sequence $S$, while $P_3$ is $3$-$\Spch$-critical exactly for packing sequences $S = (s_1, s_2, \ldots)$, where $s_1 \geq 2$. For the rest of the proof we may thus assume that $n(G) \geq 4$. 

Let $u \in V(G)$ be an arbitrary vertex of $G$. Since $G$ is a $3$-$\Spch$-critical graph, we have $\Spch(G-u)=2$ or $\Spch(G-u)=1$. The later case means that $G-u$ is a disjoint union of isolated vertices, and hence $G$ would be $2$-colorable for every packing sequence $S$. Henceforth $\Spch(G-u)=2$ holds. We now distinguish three cases with respect to the shape of $S$.

\medskip\noindent
\textbf{Case 1.} $s_1 = s_2 = 1$.\\
In this case it is clear that $\chi(G-u)=\Spch(G-u)=2$, and hence $G-u$ is a disjoint union of connected bipartite graphs and isolated vertices. If $u$ would be adjacent to vertices from at most one partition of each of the bipartite graphs, then $G$ would be $2$-colorable. Thus there exists at least one connected bipartite component of $G-u$, say $G_1$, such that $u$ has neighbors in both bipartition sets of $G_1$ and such that the subgraph induced by $V(G_1)\cup \{u\}$ contains an odd cycle $C_{2k+1}$, $k \geq 1$. Then $V(G) = V(C_{2k+1})$, for otherwise removing a vertex not belonging to the cycle would yield a $3$-colorable graph. Moreover, $E(G) = E(C_{2k+1})$, for otherwise an additional edge would yield a shorter odd cycle, so $G$ would not be $3$-$\Spch$-critical.

\medskip\noindent
\textbf{Case 2.} $s_2 \geq 2$.\\
In this case we can follow a similar line of thought than in the proof of \cite[Proposition 4.1]{klavzar-2019}. In view of Proposition~\ref{prp:2-col}, if $u$ is a vertex of a $3$-$\Spch$-critical graph $G$, then $G-u$ is a disjoint union of stars and isolated vertices. It is clear that $G-u$ must contain at least one star, say $G_1$, for otherwise $G$ would itself be a star. If $G-u$ contains more than one star, then $G$ contains $P_5$, and cannot be $3$-$\Spch$-critical, because by removing an end-vertex of $P_5$, the obtained graph would contain a $P_4$ for which $\Spch(P_4)=3$ for every packing sequence $S$ with $s_2 \geq 2$. Also, if $G-u$ has more than one isolated vertex, then removing one such vertex from $G$ yields a graph with $\Spch=3$. Thus, $G-u$ contains one star an at most one isolated vertex.

First suppose that $G-u=G_1$. Since $n(G) \geq 4$, the star $G_1$ must have at least two leaves. If $u$ is adjacent to the center of $G_1$, then since $G$ itself is not a star, $u$ is adjacent to at least one leaf of $G_1$. Removing one of the other leaves (say $v$) in $G_1$ gives a graph that contains $C_3$, and $\Spch(G-v)=3$, which is a contradiction. Therefore, $u$ is not adjacent to the center of $G_1$, hence it is adjacent to at least one leaf of $G_1$. If $S_1$ contains at least three leaves, then $G$ is not $3$-$\Spch$-critical because then we can remove a leaf and the obtained graph contains $P_4$, for which, as already noticed, $\Spch(P_4)=3$ for every packing sequence $S$ with $s_2 \geq 2$.  Hence $G_1$ must have exactly two leaves. If $u$ is adjacent to exactly one of them, we get $P_4$, and if it is adjacent to both of them, we get $C_4$. Both graphs are $3$-$\Spch$-critical for any packing sequence $S$ with $s_2 \geq 2$.

The other case to consider is when $G-u$ is a disjoint union of the star $G_1$ and an isolated vertex, say $w$. If $G_1$ has at least two leaves, then we deduce as in the subcase above that $G$ is not $3$-$\Spch$-critical. And if $G_1=K_2$, then $G$ is either $P_4$, which is $3$-$\Spch$-critical for any packing sequence $S$ with $s_2 \geq 2$, or $Z_1$  which is not $3$-$\Spch$-critical. 

\medskip\noindent
\textbf{Case 3.} $s_1 \geq 2$.\\
In this case, Proposition~\ref{prp:2-col} implies that if $u$ is a vertex of a $3$-$\Spch$-critical graph $G$, then $G-u$ is a disjoint union of copies of $K_2$ and isolated vertices. If $G-u$ has more than one $K_2$, then $G$ contains $P_5$ and cannot be $3$-$\Spch$-critical, since by removing an end-vertex of $P_5$ the obtained graph would still have a $P_4$ which is $3$-$\Spch$-colorable since $s_1 \geq 2$. Also, if $G-u$ contains isolated vertices, then $G$ contains a path $P_4$, and removing a vertex from this path yields a $P_3$ in $G-u$. Thus, $\Spch(G-u)=3$ for any packing sequence $S$ with $s_1 \geq 2$, which is a contradiction. We conclude that $G-u=K_2$ must hold and therefore $G=C_3$ or $P_3$. Both remaining graphs are $3$-$\Spch$-critical for every packing sequence with $s_1 \geq 2$.
\qed

Note that Theorem~\ref{thm:3-critical} implies that  $C_3$ is the unique graph that is $3$-$\Spch$-critical for every packing sequence $S$.

%%%%%%%%%%%%%%%%%%%%%%%%%%%%%%%%%%%%%%%%%%%%%%%%%%%%%
%%%%%%%%%%%%%%%%%%%%%%%%%%%%%%%%%%%%%%%%%%%%%%%%%%%%%
\section{On $4$-$\Spch$-critical graphs}
\label{sec:4-critical}
%%%%%%%%%%%%%%%%%%%%%%%%%%%%%%%%%%%%%%%%%%%%%%%%%%%%%
%%%%%%%%%%%%%%%%%%%%%%%%%%%%%%%%%%%%%%%%%%%%%%%%%%%%%

\begin{figure}[ht!]
\centering
 \includegraphics{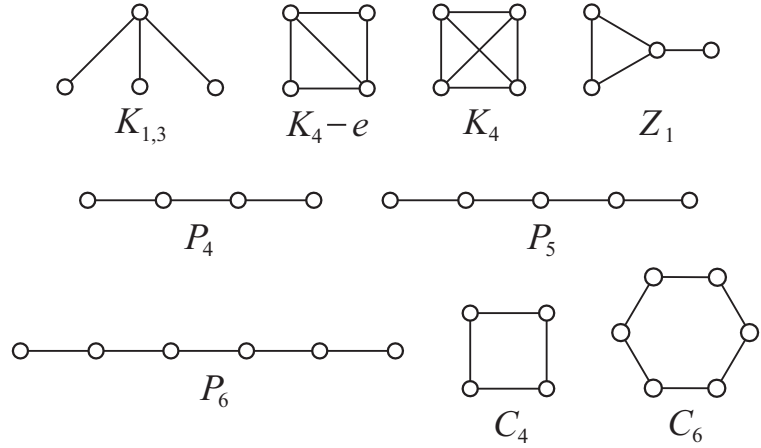}
\caption{The $4$-$\Spch$-critical graphs for packing sequences $S$ with $s_1\geq 2$}
\label{fig:4-crit}
\end{figure}

In this section we deal with the $4$-$\Spch$-critical graphs for packing sequences with $s_1\ge 2$. All critical graphs from appear in Theorem \ref{Thm-4-crit} are depicted in Fig.~\ref{fig:4-crit}.

\begin{theorem} \label{Thm-4-crit}
Let $S=(s_1,s_2,\dots)$ be a packing sequence with $s_1\geq 2$, and let $G$ be a graph.
\begin{enumerate}
\item If $s_3=2$, then $G$ is $4$-$\Spch$-critical if and only if $G \in \{K_{1,3}, C_4, Z_1, K_4-e, K_4\}$.
\item If $s_2=2$ and $s_3 \geq 3$, then $G$ is $4$-$\Spch$-critical if and only if \\ $G\in \{K_{1,3},C_4,Z_1, K_4-e, K_4, P_6,C_6\}$.
\item If $s_1=2$ and $s_2\geq 3$, then $G$ is $4$-$\Spch$-critical if and only if \\ $G\in \{K_{1,3},C_4,Z_1, K_4-e, K_4,P_5\}$.
\item If $s_1 \geq 3$, then $G$ is $4$-$\Spch$-critical if and only if $G \in \{  K_{1,3}, P_4, C_4, Z_1, K_4-e, K_4\}$.
\end{enumerate}
\end{theorem}

\proof
Let $G$ be a $4$-$\Spch$-critical graph. Then $G$ is connected by Lemma~\ref{lem:connected}. Clearly, $n(G)\geq 4$. Since $s_1\geq 2$ we have $\Delta(G)\leq 3$, for otherwise $G$ is not $4$-$\Spch$-colorable ($\Spch(G)\geq \Delta(G)+1$).

\begin{claim} \label{claim1}
Let $G$ be a $4$-$\Spch$-critical graph different from $K_{1,3}$ and let $u\in V(G)$ be such that $\Spch(G-u)=3$. If $G-u$ is disconnected, then $d_G(x)\leq 2$ for any $x\in V(G-u)$.
\end{claim}

\proof
Suppose that $G-u$ consists of at least two components. Let to the contrary $G$ contain a vertex $x\not=u$ of degree $3$, and let $G_1$ denote the component of $G$ containing $x$. Then, deleting any vertex of $G-u$ not belonging to $G_1$ yields a graph with a vertex of degree $3$, implying that $G-u$ is not $3$-$\Spch$-colorable, which is a contradiction. 
\smallqed

\noindent
\textbf{Case 1.} $s_1=2$.\\
First suppose that $\Spch(G-u)=1$ for each $u\in V(G)$. Then, clearly, $G-u$ consists of isolated vertices. Since $n(G)\geq 4$ and $d_G(u)\leq 3$, we get $K_{1,3}$. But deleting any leaf of it we get a graph which is not $1$-$\Spch$-colorable, which is a contradiction.

\smallskip
Now suppose that for each vertex $u'$ of $G$, $\Spch(G-u')\leq 2$ and there exists a vertex $u\in V(G)$ such that $\Spch(G-u)= 2$. By Proposition \ref{prp:2-col}, $G-u$ is a disjoint union of at least one copy of $K_2$ and some (possibly zero) isolated vertices. Clearly $d_G(u)\leq 3$. Since $n(G)\geq 4$, $G-u$ is disconnected and hence $d_G(u)\geq 2$. If $d_G(u)=2$, then $G-u$ consists of two components and $u$ has a neighbor in each of them, implying that $G\simeq P_k$, where $k\in \{4,5\}$. If $k=4$, then $\Spch(P_k)=3$ whenever $s_1=2$, which is a contradiction; if $k=5$, then $\Spch(P_k)=4$ when  $s_2\geq 3$, otherwise $\Spch(P_k)=3$. Thus the only critical graph for $s_2\geq 3$ is $P_5$. 

Now assume that $d_G(u)=3$. If some vertex $v$ of $G-u$ is not adjacent to $u$ in $G$, then $d_{G-v}(u)=3$, implying that $\Spch(G-v)\geq 4$ and hence $G$ is not $4$-$\Spch$-critical. If $G-u$ consists of three components, then each of these components must be $K_1$, implying that $\Spch(G-u)= 1$, which is a contradiction. Thus $G-u$ has two components and exactly one of them has two vertices. Then we get $G\simeq Z_1$, but deleting the leaf of $Z_1$ we get $C_3$ which is not $2$-$\Spch$-colorable, which is again a contradiction.

\smallskip
Finally suppose that $G$ contains a vertex $u$ such that $\Spch (G-u)=3$. Clearly, $n(G-u)\geq 3$ and $d_G(u)\leq 3$. 

Assume that $d_G(u)=3$. If some vertex $v$ of $G-u$ is not adjacent to $u$ in $G$, then $\Delta(G-v)=3$, implying that $\Spch(G-v)\geq 4$ and hence $G$ is not $4$-$\Spch$-critical.
Thus $n(G-u)=3$. If $G-u$ is disconnected, then each component of $G-u$ has at most $2$ vertices and hence $\Spch(G-u)\leq 2$, which is a contradiction. Thus $G-u$ is connected, $G-u\simeq G'\in \{P_3, C_3\}$, implying that $G\in \{ K_4-e, K_4 \}$.

Now we assume that $d_G(u)\leq 2$ and consider the following possibilities.

\medskip\noindent
\textbf{Subcase 1.1} $s_3=2$ (and also $s_1 = s_2 = 2$).\\
By Proposition \ref{prp:3-col-222}, $G-u$ consists of a disjoint union of $K_1$, $K_2$, some paths of arbitrary lengths, and of cycles of lengths divisible by $3$. 

Assume that $d_G(u)=1$. Since $G$ is connected and $n(G)\geq 4$, we infer that $G-u\simeq G'\in \{P_k:\ k\geq 3\} \cup \{C_{3k}:\ k\geq 1\}$. Connecting $u$ to an end-vertex of any path we get a path which is still $3$-$\Spch$-colorable, which is a contradiction. Connecting $u$ to a vertex of degree $2$ of any path we get $G\simeq K_{1,3}$, or a graph which is not $4$-$\Spch$-critical, since for any leaf $v\in V(G)$, such that $v$ is adjacent to a vertex of degree $2$ in $G$, we get $\Delta(G-v)=3$ implying that $\Spch(G-v)\geq 4$. Analogously, connecting $u$ to a vertex of a $C_k$, $k>3$, we get a graph which is not $4$-$\Spch$-critical. Connecting $u$ with one vertex of $C_3$ we get $G\simeq Z_1$.

Assume that $d_G(u)=2$. If $G-u$ is connected, then $n(G-u)\geq 3$. For $G-u\simeq P_3$ we get $G\simeq C_4$ or $G\simeq Z_1$, for $G-u\simeq C_3$ we get $G\simeq K_4-e$, while in any other case we get a graph which is not $4$-$\Spch$-critical or is $3$-$\Spch$-colorable (a path). Thus let $G-u$ be disconnected and consisting of two components $G_1$ and $G_2$. By Claim \ref{claim1}, $u$ is adjacent to vertices of degree $1$ only. Thus $G_1$ and $G_2$ are both paths and $u$ is adjacent to one end-vertex of $G_1$ and one endvertex of $G_2$. Then $G$ is a path, hence it is $3$-$\Spch$-colorable, which is a contradiction.

\medskip\noindent
\textbf{Subcase 1.2.} $s_2=2$ and $s_3\geq 3$.\\
By Proposition \ref{prp:3-col-223}, $n(G-u) \leq 5$. Clearly $\Delta(G-u)\leq 2$ since $\Spch(G-u)=3$ and $\Spch(G)\geq \Delta(G)+1$. Since none of $C_4$ and $C_5$ is $3$-$\Spch$-colorable, $G-u$ is a disjoint union of some copies of $K_1$, $K_2$, $P_3$, $C_3$, $P_4$, and/or $P_5$. 

Assume that $d_G(u)=1$. Since $G$ is connected and $n(G)\geq 4$, the graph $G-u$ must be one of $P_3$, $C_3$, $P_4$, and $P_5$. If $G-u \simeq C_3$, then $G\simeq Z_1$. If $u$ is adjacent to an end-vertex of $P_3$, $P_4$, or $P_5$, then we either get a $3$-$\Spch$-colorable graph or $P_6$, hence $G\simeq P_6$. If $u$ is adjacent to the central vertex of $P_3$, then $G\simeq K_{1,3}$. If $u$ is adjacent to some vertex of $P_4$ or $P_5$ of degree $2$, then we get a graph which is not $4$-$\Spch$-critical since $\Spch(G-v)=4$ for any $v\in V(G)$ with $d_G(v)=1$, and $v$ is adjacent to a vertex of degree $2$ in $G$.

Assume that $d_G(u)=2$.  If $G-u$ is not connected, then $u$ has a neighbor in two distinct components of $G-u$, implying that $\Delta(G)=2$ by Claim \ref{claim1}. Thus we get $G \simeq P_6$ since $P_4$ and $P_5$ are both $3$-$\Spch$-colorable. If $G-u$ is connected, then $G-u \in \{ P_3,C_3,P_4,P_5 \}$ since $n(G)\geq 4$. Then, for $G-u\simeq P_3$ we get $G\in \{C_4, Z_1\}$, for $G-u\simeq C_3$ we get $G\in \{K_4-e, K_4\}$, and for $G-u\simeq P_5$ we get $G\simeq C_6$; in any other case $G$ is not $4$-$\Spch$-critical.

\medskip\noindent
\textbf{Subcase 1.3.} $s_2 \geq 3$.\\
Since $s_2\geq 3$, $s_3\geq 3$ as well. Hence, by Proposition \ref{prp:3-col-223}, $n(G-u) \leq 5$. Clearly, $\Delta(G-u)\leq 2$. Since none of $P_5$, $C_4$, $C_5$ is $3$-$\Spch$-colorable, $G-u$ is a disjoint union of some copies of $K_1$, $K_2$, $P_3$, $C_3$ and/or $P_4$. 

Assume that $d_G(u)=1$. Since $n(G)\geq 4$, $G-u\simeq G'\in \{P_3,C_3,P_4\}$. If $G-u\simeq P_3$, then connecting $u$ with the central vertex of $P_3$ we get $G\simeq K_{1,3}$, otherwise connecting $u$ with an end-vertex of $P_3$ we get a $3$-$\Spch$-colorable graph $P_4$, which is a contradiction. If $G-u\simeq P_4$, then connecting $u$ with an end-vertex of $P_4$ we get $G\simeq P_5$, otherwise connecting $u$ with some vertex of degree $2$ we get a graph which is not $4$-$\Spch$-critical. And, if $G-u\simeq C_3$, we get $G\simeq Z_1$.

Assume that $d_G(u)=2$. If $G-u$ is not connected, then $u$ has a neighbor in two distinct components of $G-u$, implying that $\Delta(G)=2$ by Claim \ref{claim1}. Thus we get $G\simeq P_5$, since $P_k$ is not $4$-$\Spch$-critical for any $k\geq 6$ and $P_4$ is $3$-$\Spch$-colorable. If $G-u$ is connected, then $G-u \in \{P_3,C_3,P_4\}$ since $n(G)\geq 4$. Then, for $G-u\simeq P_3$ we get $G\in \{C_4, Z_1\}$ and for $G-u\simeq C_3$ we get $G\simeq K_4-e$; in any other case $G$ is not $4$-$\Spch$-critical.

\medskip\noindent
\textbf{Case 2.} $s_1 \geq 3$.\\
If $\Delta (G)=3$, then $n(G) = 4$, implying that $G\in \{K_{1,3}, Z_1, K_4-e, K_4\}$, for otherwise we get a graph which is not $4$-$\Spch$-colorable. If $\Delta(G)=2$, then since $G$ is connected, $G\simeq P_k$ or $C_k$, $k\geq 4$. Clearly $P_k$ is not $4$-$\Spch$-critical for any $k\geq 5$, hence we get $G\simeq P_4$. For cycles, $C_k$ is $4$-$\Spch$-colorable if and only if $k=4$, or $s_4=3$ and $k$ is divisible by $4$. And, when $k>4$, $C_k$ is not $4$-$\Spch$-critical since $\Spch(C_k-u)=4$ for any $u\in V(C_k)$ since $C_k-u$ contains a $P_4$ for which $\Spch(P_4)=4$. Thus $G\simeq C_4$. Note that the graphs $K_{1,3}$, $P_4$ and $C_4$ are all $4$-$\Spch$-critical.

\medskip\noindent
Finally note that each of the graphs from the set $\{ K_{1,3}, C_4, Z_1, K_4-e, K_4 \}$ are $4$-$\Spch$-critical for every packing sequence $S$ with $s_1\geq 2$, each of the graphs $P_6$ and $C_6$ are $4$-$\Spch$-critical for $s_2=2$ and $s_3\geq 3$, and the graph $P_5$ is $4$-$\Spch$-critical for $s_1=2$ and $s_2\geq 3$.
\qed

%%%%%%%%%%%%%%%%%%%%%%%%%%%%%%%%%%%%%%%%%%%%%%%%%%%%%
%%%%%%%%%%%%%%%%%%%%%%%%%%%%%%%%%%%%%%%%%%%%%%%%%%%%%
\section{$\Spch$-critical trees}
\label{sec:tress}
%%%%%%%%%%%%%%%%%%%%%%%%%%%%%%%%%%%%%%%%%%%%%%%%%%%%%
%%%%%%%%%%%%%%%%%%%%%%%%%%%%%%%%%%%%%%%%%%%%%%%%%%%%%

Let $S$ be a packing sequence with $s_1 = s_2 = 1$. Then every bipartite graph, in particular every tree, admits a $2$-$\Spch$-coloring. It follows that $K_2$ is the only $\Spch$-critical graph for such a packing sequence. On the other hand, if $s_2\ge 2$, then the situation is more interesting already on trees  as the next result which extends \cite[Proposition~5.1]{klavzar-2019} asserts.

\begin{proposition}
If $k\ge 2$ and $S$ is a packing sequence with $s_2\ge 2$, then there exists a $k$-$\Spch$-critical tree. 
\end{proposition}

\proof
Suppose first that $s_1 \ge 2$ (and, of course, $s_2\ge s_1$). Then $K_{1,k-1}$ is a required $k$-$\Spch$-critical tree.

Assume in the rest that $s_1 = 1$ and (and $s_2\ge 2$). Let $T_k$ be the tree obtained from $K_{1,k-1}$ with the central vertex $u$ and leaves $w_1, \ldots, w_{k-1}$, by attaching $k-2$ leaves to each of the vertices $w_i$, $i\in [k-1]$. In particular, $T_2 = K_2$ and $T_3 = P_5$. 

We claim that $\Spch(T_k) = k$. Let $c$ be an arbitrary $\Spch$-coloring of $T_k$. If $c(w_i) = 1$ holds for some $i\in [k-1]$, then the $k-1$ neighbors of $w_i$ must receive pairwise different colors, hence $c$ uses at least $k$ colors. On the other hand, if $c(w_i)\ne 1$ for each $i\in [k-1]$, then $c$ uses $k-1$ colors on the vertices $w_i$ and hence at least $k$ colors all together. This shows that $\Spch(T_k) \ge k$. On the other hand, setting $c(w_i) = i+1$, $i\in [k-1]$, and coloring every other vertex with color $1$ yields $\Spch(T_k) \le k$. This proves the claim. 

We have thus seen that $\Spch(T_k) = k$. If $T_k$ is $k$-$\Spch$-critical, then we are done. Otherwise, using Lemma~\ref{lem:leaf}, remove leaves of $T_k$ one by one until a $k$-$\Spch$-critical tree is obtained.         
\qed

In \cite{klavzar-2019}, $k$-$\Spch$-critical caterpillars were investigated for $S=(1,2,3, \dots)$. Here we focus on existence of $k$-$\Spch$-critical caterpillars for some packing sequences with $s_1=1$. Recall that, for $s_2=1$, the only $\Spch$-critical graph is $K_2$, thus we consider $s_2\geq 2$.

\begin{proposition}
Let $k$ be a positive integer and $S=(1,s_2^{k-1})$ a packing sequence with $s_2\geq 2$. Then a $k$-$\Spch$-critical caterpillar exists if and only if $k\leq s_2+2$.
\end{proposition}

\proof
Let $T$ be a caterpillar. Since any path has a $(s_2^{s_2+1})$-packing coloring repeating the coloring pattern $2,3,\dots, s_2+2$, we can color vertices of the spine of $T$ with colors $2,3,\dots s_2+2$ and then color all the leaves of $T$ with color $1$. Thus $\Spch(T)\leq s_2+2$ for an arbitrary caterpillar $T$.

On the other hand, any path of length greater than $s_2$ has no $(s_2^{s_2})$-packing coloring. Then, considering any path $P$ of length greater than $s_2$ and attaching at least $s_2+2$ leaves to each vertex of $P$ we get a caterpillar $T$ with $\Spch(T)=s_2+2$. Iteratively applying Lemma \ref{lem:leaf} to the leaves of $T$ we find an $(s_2+2)$-$\Spch$-critical caterpillar. Continuing deleting leaves of $T$ in this manner we can also get a $k$-$\Spch$-critical caterpillar for any $k\leq s_2+1$.
\qed

Now we construct an explicit $k$-$\Spch$-critical caterpillars for $S=(1,s_2^{k-1})$ and any $k\leq s_2+2$, $k\geq 2$.

\begin{figure}[ht!]
\centering
 \includegraphics{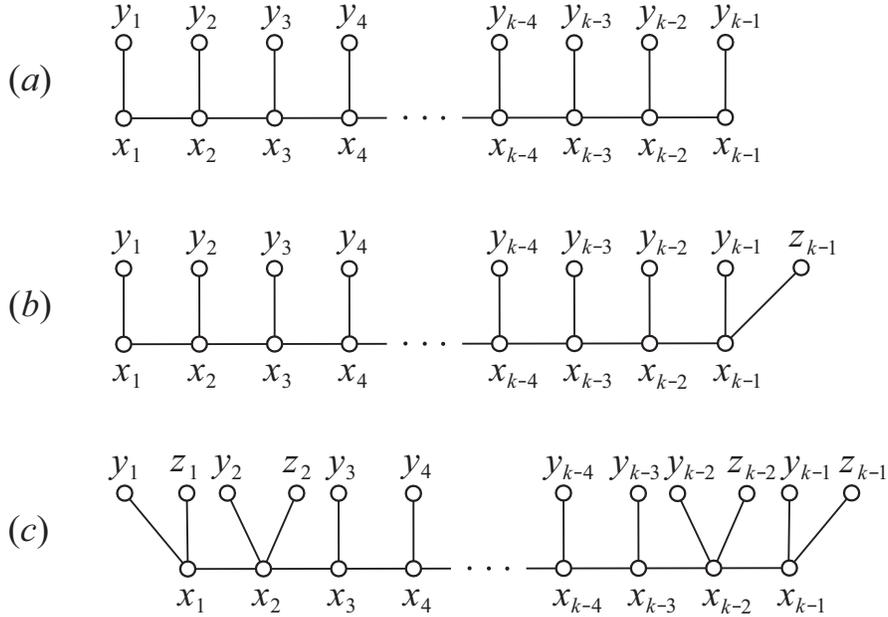}
\caption{Caterpillars}
\label{fig:caterpillars}
\end{figure}

\begin{example} \rm
Let $S=(1,s_2^{k-1})$ be such that $k\leq s_2$. Let $G_1$ be a caterpillar consisting of a spine $P$ of length $k-2$ and adding one leaf to each vertex of $P$ (see Fig. \ref{fig:caterpillars} $(a)$). We show that $G_1$ is $k$-$\Spch$-critical. First, since each color different from $1$ can be used for only one vertex of $G_1$ and color $1$ can be used on at most $k-1$ vertices, we have $\Spch(G_1)\geq k$. On the other hand, we can color vertices of $P$ with colors $2,3,\dots, k$ and all leaves of $G_1$ with color $1$, implying that $\Spch(G_1)=k$. 

Now we show that $\Spch(G_1-x)<k$ for any $x\in V(G_1)$. Deleting any leaf $x$ of $G_1$ we get a graph in which we color all leaves of $G_1-x$ and the neighbor of $x$ in $G_1$ with color $1$, and we color the remaining $k-2$ vertices of $P$ with mutually distinct colors $2,3,\dots, k-2$, implying that $\Spch(G_1-x)<k$. If $x\in V(P)$, then $G_1-x$ is disconnected and consisting of one isolated vertex and one or two caterpillars $C_1,C_2$. Clearly, each of $C_1,C_2$ has a spine of length smaller than $k-1$, implying that $\Spch(C_i)<k$ for $i=1,2$. Therefore $G_1$ is $k$-$\Spch$-critical.
\end{example}

\begin{example} \rm
Let $S=(1,s_2^{k-1})$ be such that $k=s_2+1$. Let $G_2$ be a caterpillar consisting of a spine $P=x_1,x_2,\dots, x_{k-1}$, adding one leaf $y_i$ to $x_i$ for each $i=1, \dots ,k-1$ and adding one more leaf $z_{k-1}$ to $x_{k-1}$ (see Fig. \ref{fig:caterpillars} $(b)$). We show that $G_2$ is $k$-$\Spch$-critical. First, since color $1$ can be used for an independent set of $G_2$, we have to color either $x_i$ or all leaves adjacent to $x_i$ with some of the colors $2,3,\dots, k$ for each $i=1, \dots, k-1$. And, since we can use only one such color twice (for $y_1$, and $y_{k-1}$ or $z_{k-1}$), we have $\Spch(G_2)\geq k$. On the other hand, coloring vertices of $P$ with colors $2,3,\dots, k$ and coloring all leaves of $G_2$ with color $1$, we get an $k$-$\Spch$-coloring of $G_2$, implying that $\Spch(G_2)=k$.

Now we show that $\Spch(G_2-x)<k$ for any $x\in V(G_2)$. Deleting any leaf not adjacent to $x_{k-1}$ we get a graph in which we color all leaves of $G_2-x$ and the neighbor of $x$ in $G_2$ with color $1$, and we color the remaining $k-2$ vertices of $P$ with mutually distinct colors $2,3,\dots, k-1$, implying that $\Spch(G_2-x)<k$. If $x$ is adjacent to $x_{k-1}$, say, $x=z_{k-1}$, we color $y_1$ and $y_{k-1}$ with color $2$, all internal vertices of $P$ (one by one) with colors $3,4,\dots ,k-1$ and all the remaining vertices with color $1$. Thus, in this case, $\Spch(G-x)<k$. If $x\in V(P)$, then $G_2-x$ is disconnected and it consists of one or two isolated vertices and one or two caterpillars $C_1,C_2$. Clearly, each of $C_1,C_2$ has a spine of length smaller than $k-1$, implying that $\Spch(C_i)<k$ for $i=1,2$. Therefore $G_1$ is $k$-$\Spch$-critical. 
\end{example}

\begin{example} \rm
Let $S=(1,s_2^{k-1})$ be such that $k=s_2+2$. Let $G_3$ be a caterpillar consisting of a spine $P=x_1,x_2,\dots, x_{k-1}$, adding one leaf $y_i$ to $x_i$ for each $i=3, \dots ,k-3$, and adding two leaves $y_i$, $z_i$ to $x_i$ for each $i=1,2,k-2,k-1$ (see Fig. \ref{fig:caterpillars} $(c)$). We show that $G_3$ is $k$-$\Spch$-critical. First, since color $1$ can be used for an independent set of $G_3$, we have to color either $x_i$ or all leaves adjacent to $x_i$ with some of the colors $2,3,\dots, k-1$ for each $i=1, \dots, k-1$. It is straightforward to check that we always need at least $k-1$ such colors, implying that $\Spch(G_3)\geq k$. On the other hand, coloring vertices of $P$ with colors $2,3,\dots, k$ and coloring all leaves of $G_3$ with color $1$, we get an $k$-$\Spch$-coloring of $G_3$, implying that $\Spch(G_3)=k$.

Now we show that $\Spch(G_3-x)<k$ for any $x\in V(G_3)$. Deleting any leaf $x=y_i$ for some $i=3,4, \dots, k-3$, we get a graph in which we can color all leaves of $G_3-x$ and the neighbor of $x$ in $G_3$ with color $1$, and we color the remaining $k-2$ vertices of $P$ with mutually distinct colors $2,3,\dots, k-1$, implying that $\Spch(G_3-x)<k$. If, up to symmetry, $x=y_1$, then we can color  $x_1$ and all leaves of $G_3-x$ different from $z_1$ with color $1$, $z_1$ and $x_{k-1}$ with color  $2$, and the remaining $k-3$ vertices of $P$ one by one with colors $3, \dots, k-1$, implying that $\Spch(G_3-x)<k$. Analogously, if, up to symmetry, $x=y_2$, then we can color $x_2$, $x_{k-1}$ and all leaves of $G_3-x$ different from $z_2$  with color $1$, $x_1$ and $y_{k-1}$ with color $2$, $z_2$ and $z_{k-1}$ with color $3$, and the remaining $k-4$ vertices of $P$ one by one with colors $4,5,\dots, k-1$, implying that $\Spch(G_3-x)<k$. If $x\in V(P)$, then $G_3-x$ is disconnected and it consists of one or two isolated vertices and one or two caterpillars $C_1,C_2$. Clearly, each of $C_1,C_2$ has a spine of length smaller than $k-1$, implying that $\Spch(C_i)<k$ for $i=1,2$. Therefore $G_3$ is $k$-$\Spch$-critical. 
\end{example}

\section*{Concluding remarks}

The criticallity studied in~\cite{klavzar-2019} for the packing chromatic number and the criticallity investigated in this paper for the $S$-packing chromatic number refer to vertex-deleted subgraphs. An equally legal criticallity concept is the one with respect to arbitrary subgraphs, equivalently with respect to edge-deleted subgraphs. The seminal study~\cite{bresar-2019+} on the latter concept for the packing chromatic number  has been done independently and at about the same time as~\cite{klavzar-2019}. It would hence be natural to study also the edge-deleted critical graphs in the general context of $S$-packing colorings. 

It was stated in~\cite{klavzar-2019} that it would be interesting to classify vertex-transitive, $\pch$-critical graphs. Here we extend this claim by stating that it would also be of interest to classify vertex-transitive, $\Spch$-critical graphs for each of the packing sequence $S$.

\section*{Acknowledgements}

We acknowledge the financial support from the project GA20-09525S of the Czech Science Foundation (first author) and from the Slovenian Research Agency (research core funding No.\ P1-0297 and projects J1-9109, J1-1693, N1-0095) (second and third author). 

%%%%%%%%%%%%%%%%%%%%%%%%%%%%%%%%%%%%%%%%%%%%%%%%%%%%%%%%%%%%%%%%%%%%%
%%%%%%%%%%%%%%%%%%%%%%%%%%%%%%%%%%%%%%%%%%%%%%%%%%%%%%%%%%%%%%%%%%%%%

\end{document}